\documentclass[11pt]{amsart}

\usepackage{amssymb}
\usepackage{latexsym}
\input xypic
\xyoption {all}

\usepackage{color}

\theoremstyle{plain}
\newtheorem{thm}{Theorem}[section]

\newtheorem{prop}[thm]{Proposition}

\newtheorem{lemma}[thm]{Lemma}

\newtheorem{prob}[thm]{Problem}

\newtheorem{ex}[thm]{Example}

\theoremstyle{definition}
\newtheorem*{ack}{Acknowledgements}
\newtheorem{rema}[thm]{Remark}
\newtheorem{defn}[thm]{Definition}

\numberwithin{equation}{section}

\newcommand{\de}{\delta}

\newcommand{\e}{\varepsilon}

\newcommand{\bbN}{\mathbb{N}}

\newcommand{\supp}{\operatorname{supp}}

\newcommand{\disp}{\displaystyle}
\newcommand{\lb}{\label}

\newcommand{\wtw}{if and only if }

\begin{document}

\dedicatory{Dedicated to the memory of James E. Jamison}

\title{On sign embeddings and narrow operators on $L_2$}


\author{Beata~Randrianantoanina}

\address{Department of Mathematics\\
Miami University\\
Oxford, OH 45056, USA}

\email{\texttt{randrib@miamioh.edu}}




\begin{abstract}
The goal of this note is two-fold. First we present a brief overview of ``weak'' embeddings, with a special emphasis on sign embeddings which were introduced by H. P. Rosenthal in the early 1980s. We also discuss the related notion of narrow operators, which was introduced by A. Plichko and M. Popov in 1990. We give examples of applications of these notions in the geometry of Banach spaces and in other areas of analysis. We also present some open problems. 

In the second part  we prove that  Rosenthal's celebrated characterization of narrow operators  on $L_1$ is also true for operators on $L_2$. This answers, for $p=2$, a question posed by  Plichko and  Popov in 1990. For $1<p<2$ the problem remains open.
\end{abstract}

\maketitle

\section{A brief overview of ``weak'' embeddings and narrow operators}
\label{intro}

Isomorphisms and isomorphic embeddings are a ``best'' kind of an operator (in this note, an operator will always mean a bounded linear operator) between function spaces and Banach spaces in the sense that they preserve many properties of the spaces on which they operate. A natural question to consider is whether it is possible to identify a class (or classes) of operators that is wider than isomorphic embeddings, and still preserves important properties of spaces on which they operate, or whose existence would imply an existence of an isomorphic embedding of a substantial subspace of the domain. It turns out that such classes do exist, at least for some sufficiently good Banach or function spaces. They are sometimes referred to by an umbrella term of ``weak embeddings''. They became an interest of many mathematicians since late 1970s.

Throughout this note all Banach spaces can be  real or complex. We use standard notation, we refer the reader to \cite{AK} for all undefined notions and notation.

\subsection{Definitions of various notions of ``weak'' embeddings.}

The first notion of this type, a {\it semi-embedding} (an injective operator  is a semi-embedding if the image of the closed unit ball of the domain space is closed) was introduced by Lotz, Peck and Porta
\cite{LPP79}, who used it to study properties of $c_0$ and $C(K)$, for $K$ compact. In particular, they showed that a compact space $K$ is scattered if and only if all one-to-one Tauberian operators from $C(K)$ into arbitrary Banach spaces are isomorphisms, which answered a question of Kalton and Wilansky \cite{KW76}. They also
 proved that every semi-embedding of $C[0, 1]$ is an isomorphism when restriced to a suitable
complemented subspace isomorphic to $C[0,1]$.

In the early 1980s two additional notions of weak embeddings were introduced and studied.
Bourgain and Rosenthal \cite{BR83} defined {\it $G_\de$-embeddings} as injective operators such that the image of every closed bounded set is a $G_\de$-set, and  Rosenthal \cite{Ros81/82} defined {\it sign embeddings} as injective operators such that there exists $\de>0$ with $\|Tx\|\ge \de\|x\|$ for every sign $x$ in $X$, where, if $X$ is a K\"othe function space on $[0,1]$, e.g. $X=L_p[0,1]$ for $1\le p\le\infty$, a sign $x$ is a measurable function in $X$ which only attains values from the set $\{1,-1,0\}$, i.e. $x=\bold{1}_{A_1}-\bold{1}_{A_2}$ for some measurable disjoint sets $A_1$ and $A_2$ (we use the symbol $\bold{1}_{A}$ to denote the characteristic function of the set $A$). 

\subsection{Examples of applications of weak embeddings.}

The study of weak embeddings has many applications in the geometry of Banach spaces and in related areas. We do not have the space here to list all such applications in the literature, but we do want to mention some examples.

As one of the first applications of sign embeddings Bourgain and Rosenthal \cite{BR83} deduced a stronger version of the
classical theorem of Menchoff \cite{M1916} in harmonic analysis that there exists a singular probability
measure on the circle with Fourier coefficients vanishing at
infinity, and the measure may be chosen to be singular with respect to any pre-assigned probability measure with non-vanishing Fourier coefficients.

As another application, Bourgain and Rosenthal \cite{BR83} also proved that if $X$ is a subspace of $L_1$, then either $L_1$ embeds in $X$, or $\ell_1$ embeds in $L_1/X$. Note that in some sense this is the optimal result, since later Talagrand \cite{T90} showed an example of a subspace $Z$ of $L_1$ so that $L_1$ does not embed in neither $Z$ nor $L_1/Z$.

The study of sign embeddings on $L_1$ allowed Rosenthal \cite{Ros84} to prove the following analog of  James's   fundamental $\ell_1$ Distortion Theorem.
\begin{thm}\lb{Ros1}
Let $X$ be a subspace of $L_1$, such that $X$ is isomorphic to $L_1$. Then, for every $\e>0$, there exists a subspace $Y$ of $X$, so that $Y$ is $(1+\e)$-isomorphic to $L_1$.
\end{thm}

Theorem~\ref{Ros1} is somewhat surprising and quite remarkable, since Lindenstrauss and Pe\l czy\'nski \cite{LP71} proved that the space $L_1$ is distortable, that is there exists an equivalent renorming $|||\cdot|||$ of $L_1$ and a $\de>0$ so that no subspace of
$X=(L_1, |||\cdot|||)$ is $(1+\de)$-isomorphic to $L_1$ with the usual norm.

More recently, Dosev, Johnson and Schechtman \cite{DJS} proved another result about the structure of sign embeddings on $L_p$ for $1<p<2$ (see Theorem~\ref{DJS} below), and used it to characterize commutators on $L_p$ for $1\le p<\infty$.

Many other applications of weak embeddings have been obtained, in particular in connection with the Radon-Nikodym Property, the Krein-Milman Property, and the Daugavet Property.

\subsection{Definition of narrow operators}

In 1990, Plichko and Popov \cite{PP} coined the term {\it narrow operator} to mean an operator which is, in a sense that we explain below, an ``opposite'' of a  sign embedding.

Before we give the formal definition of narrow operators, let us clarify a somewhat confusing usage of the term {\it sign embedding} in the literature. As we mentioned above, sign embeddings were first named and defined by Rosenthal in \cite{Ros81/82} as injective operators on $L_1$ (the definition also makes sense for  $L_p$, $1\le p\le\infty$) such that there exists $\de>0$ with $\|Tx\|\ge \de\|x\|$ for every sign $x$ in $L_1$. In \cite{MPRS} the term sign embedding was used for operators on $L_p$, $1\le p\le\infty$, which satisfy the same condition as sign embeddings of \cite{Ros81/82}, but are not necessarily injective. In \cite{Ros84} Rosenthal defined a similar notion which he called a {\it sign-preserving operator} (or a {\it norm-sign-preserving operator} in \cite{GR83}). An operator $T$, which is not necessarily injective, on  $L_1$ (or, analogously, on $L_p$, $1\le p\le\infty$) is called {\it sign-preserving} if there exists $\de>0$  and there exists a set $A$ of positive measure so that $\|Tx\|\ge \de$ for every sign $x$ with $\supp x=A$ and $\int x \ d\mu=0$. However, in \cite{DJS} the operators satisfying the latter condition were called sign embeddings instead of sign-preserving. The reason that this does not lead to too much confusion, despite the fact that the above three notions are distinct from each other, is that Rosenthal in \cite{Ros83} proved that if there exists a sign-preserving operator $T:L_1\to X$, then there exists a set $B$ of positive measure, so that $T\big|_{L_1(B)}$ is a sign embedding in the sense of \cite{Ros81/82}. An analogous result is valid in $L_p$, $1\le p<\infty$, see  \cite{MP06,PRbook} for a full discussion, examples illustrating differences between the different definitions, and proofs. 
Thus, since for every set $B\subseteq [0,1]$ of positive measure, $L_p[0,1]$ and $L_p(B)$ (with the measure scaled so that $B$ is of measure 1) are isometrically isomorphic to each other, all results about properties and existence of sign embeddings on $L_p$ are valid in the sense of each of the above definitions.

In the sequel we will follow the convention of \cite{MPRS}.

\begin{defn} \label{se}
Let $1\le p<\infty$, and $X$ be a Banach space. We say that a (not necessarily injective) operator $T:  L_p\to X $ is a \textit{sign embedding} if there exists $\de>0$ so that $\|Tx\|\ge \de\|x\|$ for every sign $x$ in $L_p$, where a sign is any measurable function on $[0,1]$ with values in the set $\{1,-1,0\}$.
\end{defn}

Plichko and Popov \cite{PP} defined {\it narrow operators} as the operators which are not sign-preserving, that is they introduced the following definition.

\begin{defn} \label{dnarrow} \cite{PP}
Let $1\le p<\infty$, and $X$ be a Banach space. We say that an operator $T:  L_p\to X $ is  \textit{narrow} if  for every $\varepsilon>0$ and every measurable set $A\subseteq [0,1]$ there exists  a sign $x\in L_p$ with $\supp x =A$ and $\int x \, d \mu = 0$, so that $\|Tx\|<\varepsilon$.  
\end{defn}

Before the work of Plichko and Popov, the study of sign embeddings was more a study of ``sign embeddability'' and mainly concentrated on questions of the type: if there exists a sign embedding of $X$ into $Y$, does there exist an isomorphic embedding of $X$ into $Y$? Plichko and Popov shifted the focus to the study of properties of the operators themselves. 

It is easy to see that every compact operator is narrow, cf. \cite{PP,PRbook}.

 On the other hand if $T:L_p[0,1]\to X$, $1\le p<\infty$, is a narrow operator, then for each subset $A$ of $[0,1]$ of positive measure and for each $\e>0$, there exists a subspace $E$ of $L_p(A)$ so that $E$ is isomorphic to $L_p$ and 
such that the restriction $T_1=T\big|_E$  of $T$ to the subspace $E$ is a compact operator and $\|T_1\|<\e$
(for $p=1$ see \cite{Ros84}, for general $p$ see \cite[Proposition~2.19]{PRbook}).

However  in general  the class of narrow operators is much larger than that of compact operators. All  spaces  $X$  that admit a non-compact operator from $L_p$ to $X$ also admit a narrow non-compact operator from $L_p$ to $X$ \cite[Corollary~4.19]{PRbook}. Actually, even a stronger fact is true. 
For all $1\le p\le\infty$, there exists a non-strictly singular operator $T:L_p\to L_p$, so that $T$ is narrow.
Indeed, the conditional expectation from $L_p([0,1]^2)$ onto the subspace of all functions depending on the first variable only is such an operator,  \cite{PP,K09}, see also \cite[Section~4.2]{PRbook}. 

\subsection{Characterizations of narrow operators.}

It is a natural and a very interesting question to find conditions on an operator $T$ which would characterize that $T$ is  narrow.

The most complete answer to this question was proved by Rosenthal \cite{Ros84} for operators from $L_1$ to $L_1$.

\begin{thm}  \label{Ros}
An operator $T:L_1\to L_1$ is narrow if and only if for each measurable set $A \subseteq [0,1]$ the restriction $T \bigl|_{L_1(A)}$ is not an isomorphic embedding.
\end{thm}

We note that in Theorem~\ref{Ros} it is important that the codomain space is $L_1$. Indeed, it follows from \cite{T90} that there exists a Banach space $X$ and a non-narrow operator $T: L_1\to X$ so that for each measurable set $A \subseteq [0,1]$ the restriction $T \bigl|_{L_1(A)}$ is not an isomorphic embedding (cf.  \cite[Corollary~8.25]{PRbook}). For the general codomain space Bourgain and Rosenthal \cite{BR83} proved the following result.
\begin{thm}  \label{BR}
Let $X$ be any Banach space. Then every $\ell_1$-strictly singular operator $T:L_1\to X$ is narrow. 
\end{thm}
Recall, that if $E, F, Z$ are infinite-dimensional Banach spaces then an operator $T:E\to F$ is called {\it $Z$-strictly singular} (resp. {\it strictly singular}) if for every subspace $E_1\subseteq E$ so that $E_1$ is isomorphic to $Z$ (resp. so that $E_1$ is infinite-dimensional), the restriction $T\big|_{E_1}$ is not an isomorphic embedding.  

Motivated by Theorems~\ref{Ros} and \ref{BR}, it is natural to ask the following questions, which were originally stated in \cite{PP}, see also \cite{PRbook} for a more detailed discussion.

\begin{prob} \label{pr:gen}
Suppose $1 \le p <\infty$,  $X$ is a Banach space, and an operator $T:L_p\to X$ satisfies one of the following (progressively stronger) conditions.

(a) $L_p$-strictly singular, or

(b) $\ell_2$-strictly singular, or

(c) strictly singular.

 Does it follow that $T$ is narrow?
\end{prob}

We start from a simple observation that for operators from $L_2$ to $L_2$ the answer to Problem~\ref{pr:gen}  (when $p=2$, conditions (a),(b),(c) are all the same) is positive since
 all strictly singular operators from $L_2$ to $L_2$ are compact and, therefore, narrow.

 When  $p > 2$ and $X=L_p$,  the answer to Problem~\ref{pr:gen}(a) is negative, as illuminated in the following  example.

 \begin{ex} \lb{Sch}
Let $p>2$ and  $T = S \circ J$ where $J: L_p \to L_2$ is the inclusion embedding and $S: L_2 \to L_p$ is an isomorphic embedding. Then $T$ is $L_p$-strictly singular and not narrow. (Also, for all $p>r\ge 1$, the inclusion $J:L_p\to L_r$ is not narrow.)
 \end{ex}

When $p=1$ and $X=L_1$, the answer to Problem~\ref{pr:gen}(a) is positive and follows from a result of Enflo and Starbird \cite{ES79}, cf. also  Rosenthal's Theorems~\ref{Ros1} and \ref{Ros}  stated above, which were proved in \cite{Ros84}. As we mentioned earlier, Talagrand \cite{T90} gave an example of a Banach space $X$ for which the answer to Problem~\ref{pr:gen}(a) is negative when $p=1$.

When $1< p< 2$, and $X=L_p$,  the answer to Problem~\ref{pr:gen}(a) is positive
and was first explicetly stated  by Bourgain \cite[Theorem~4.12, item 2]{Bou81}, who remarked that it can be deduced from the proof of a related result of Johnson, Maurey, Schechtman and Tzafriri \cite{JMST}, even though it is not explicitely stated there (details of  necessary adjustments were not presented).
 Recently Dosev, Johnson and Schechtman \cite{DJS} proved the following strengthening of the related result from  \cite{JMST}, which gives a very strong positive answer to Problem~\ref{pr:gen}(a) for $1< p< 2$, and $X=L_p$.

\begin{thm} \label{DJS}
For each $1 < p < 2$ there is a constant $K_p$ such that if the operator $T:L_p\to L_p$ is  not narrow  (and in particular, if $T$ is an isomorphism) then there is a $K_p$-complemented subspace $X$ of $L_p$ which is $K_p$-isomorphic to $L_p$ and such that $T\big|_X$ is a  $K_p$-isomorphism and $T(X)$ is $K_p$ complemented in $L_p$.
\end{thm}

It follows from Example~\ref{Sch}, that for $1<p<2$ and a general Banach space $X$, the answer to  Problem~\ref{pr:gen}(a) can be negative. We are not aware of any work in the literature characterizing which Banach spaces would yield an affirmative answer for $1<p<2$.

 Problem~\ref{pr:gen}(b) is open in general, but there are two (incomparable) situations when the answer is known to be positive.

 Flores and Ruiz \cite{FR03} proved the following result.

\begin{thm} \label{FR}
Let $1\le p<\infty$, and $F$ be an order continuous Banach lattice. Then every regular $\ell_2$-strictly singular operator $T: L_p\to F$ is narrow.
\end{thm}

Recall that an operator between Banach lattices is called {\it regular} if it is a difference of two positive operators.

For the codomain space which is not necessarily a Banach lattice, Mykhaylyuk,  Popov,   Schechtman, and the present author \cite{MPRS} proved the following result.

\begin{thm} \label{MPRS}
Let  $1 < p < \infty$,   and $X$ be a  Banach space  with an unconditional basis. Then  every  $\ell_2$-strictly singular operator $T :L_p\to X$ is narrow.
\end{thm}

We are not aware of any results in the literature  concerning Problem~\ref{pr:gen}(c) beyond what we mentioned above. We do consider it a very interesting problem.

Another  interesting problem, which was first posed by Plichko and Popov in \cite{PP}, is whether  an analogue of Theorem~\ref{Ros}  is valid for $1 < p \leq 2$.

\begin{prob} \label{pr:1}
Suppose $1 < p \leq 2$, and an operator $T:L_p\to L_p$ is such that for every  measurable set $A \subseteq[0,1]$ the restriction $T \bigl|_{L_p(A)}$ is not an isomorphic embedding. Does it follow that $T$ is narrow?
\end{prob}

The remainder of this paper is devoted to the positive answer to Problem~\ref{pr:1} for the case when $p=2$.  

We hope that this brief overview will encourage the readers to answer some of the many remaining open problems in this area of study. For example, it is not known whether every narrow operators has an invariant subspace. It also may be interesting to investigate whether semi-embeddings  and/or $G_\de$-embeddings may be useful in the setting of spaces of analytic functions, and whether there exists a meaningful analog of sign embeddings and narrow operators in this setting; we note that in \cite{MPRS} narrow operators on $L_p$, $1<p<2$, were characterized using functions with a controlled growth, instead of  sign functions. 
We  refer the reader to \cite{PRbook} where many more open problems and directions of study are discussed.

\section{New results}
\label{new}

In this section we prove that Problem~\ref{pr:1} has a positive answer when $p=2$, i.e. we prove the following result.

\begin{thm}\lb{L2main}
An operator $T:L_2\to L_2$ is narrow if and only if for each measurable set $A \subseteq [0,1]$ the restriction $T \bigl|_{L_2(A)}$ is not an isomorphic embedding.
\end{thm}

Before we start the proof of Theorem~\ref{L2main} we  present a slight weakening of the notion of a narrow operator and a  structural result proved in \cite{MPRS}, both of which will be useful for our proof.

\begin{defn} \label{d:somewhatnar}
Let $1\le p<\infty$ and $X$ be a Banach space. We say that an operator $T: L_p\to X$  is  \textit{somewhat narrow} if for each $A\subseteq [0,1]$ of positive measure,  and each $\varepsilon > 0$,  there exists a set $B \subseteq A$ of positive measure and a sign $x$ with $\supp x= B$ such that $\|Tx\| < \varepsilon \|x\|$.
\end{defn}

Obviously, each narrow operator is somewhat narrow. The inclusion embedding $J: L_p \to L_r$ where $1 \leq r < p < \infty$ is an example of a somewhat narrow operator which is not narrow.
However, for operators from $L_p$ to $L_p$, where $1\le p\le 2$,  these two notions are equivalent as was shown in \cite[Theorem~2.2]{MPRS}.

\begin{thm} \label{A}
Let $1 \leq p \leq 2$. Then an operator $T:L_p\to L_p$ is  somewhat narrow \wtw it is narrow.
\end{thm}

We observe that the following fact follows directly from our defintions.

\begin{lemma}\lb{signsomewhatnarrow}
Let $1\le p<\infty$ and $X$ be a Banach space. An operator $T:L_p\to X$ is somewhat narrow if and only if for each set $A$ of positive measure, the   restriction $T \bigl|_{L_p(A)}$ is not a sign embedding in the sense of Definition~\ref{se}.
\end{lemma}

We denote by
$(\overline{h}_n)_{n=1}^\infty$, the $L_\infty$-normalized Haar system, i.e. for $n=0,1,2,\dots$
and $k=1,\dots,2^n$,
$$\overline{h}_{2^n+k} = {\bold{1}}_{I_{n+1,2k-1}}-{\bold{1}}_{I_{n+1,2k}},$$
where $I_{n,k}$ denotes the dyadic interval ${\disp I_{n,k}=\Big[\frac{k-1}{2^n}, \frac{k}{2^n}\Big)}.$ Note that in this notation, $\overline{h}_{1}={\bold{1}}_{[0,1)}$ is the constant-1 function on $[0,1)$.

We denote by
 $(h_n)$ and $(h_n^*)$, the $L_p$- and $L_q$-normalized Haar functions respectively,
  where $1 \leq p< \infty$ and $1/p + 1/q = 1$.

We are now ready to state the main ``structural'' result which was proved in \cite[Proposition~3.1]{MPRS}.

\begin{prop}\label{structure} \cite{MPRS}
Suppose $1 \leq p< \infty$, $X$ is a Banach space with a basis $(x_n)$, $T: L_p\to X$ is such that there exists $\de>0$ so that $\|Tx\|\geq 2 \delta$ for each  sign $x \in L_p$
with $ \supp x=[0,1]$ and $\int x\ d\mu=0$. Then for each $\varepsilon>0$ there exist an operator $S : L_p\to X$, a normalized block basis $(u_n)$ of $(x_n)$ and real numbers $(a_n)$ such that
\begin{enumerate}
  \item $S h_n = a_n u_n$ for each $n \in \mathbb N$ with $a_1 = 0$;
  \item $\|Sx\|\geq \delta$, for each  sign $x \in L_p$
with $ \supp x=[0,1]$ and $\int x\ d\mu=0$;
\item there exists a measure preserving homeomorphism $\sigma:[0,1]\to[0,1]$, so that $\|Sx\|\leq 
\|T(x\circ \sigma)\| + \varepsilon$ for every $x \in L_p$ with $\|x\|=1$;
    \item there are finite codimensional subspaces $X_n$ of $L_p$ such that $\|Sx\|\leq \|T(x\circ \sigma)\| + 1/n$ for every $x \in X_n$ with $\|x\|=1$.
\end{enumerate}
If, moreover, $\|Tx\|\geq 2 \delta \|x\|$ for every sign $x$, then
$|a_n|\geq\delta$  for each $n \geq 2$.
\end{prop}

\begin{rema}
We note that the statement of item (3) in \cite[Proposition~3.1]{MPRS} was formulated slightly differently
than  above. In \cite{MPRS} we said that there exists a sign
preserving isometry $V$ of $L_p[0,1]$   so that $\|Sx\|\leq 
\|T(Vx)\| + \varepsilon$ for every $x \in L_p$ with $\|x\|=1$. 
It is well known (and easy to verify) that every sign preserving isometry $V$ of $L_p[0,1]$,    for $1 \leq p< \infty$, is of the form $Vx=\omega\cdot x\circ\sigma$, where $\sigma$ is a 
measure preserving homeomorphism of $[0,1]$ onto itself, and $\omega$ is a weight function whose modulus is equal to 1 a.e. 

In fact, the proof in \cite{MPRS}
that  an isometry $V$ in item (3) exists, is really a proof that there exists an appropriate
measure preserving homeomorphism $\sigma$ of $[0,1]$ onto itself, and defining
$V(x)=x\circ \sigma$, for $x\in L_p[0,1]$ (using the weight constantly equal to 1).
In the sequel it will be convenient to use the explicit form of the isometry $V$, so we elected to include it  in the statement of Proposition~\ref{structure}.
 \end{rema}

We are now ready for the main ingredient of the proof of Theorem~\ref{L2main}.

\begin{prop}\lb{2notsignpr}
Suppose $T:L_2\to L_2$ is such that for every  measurable set $A \subseteq[0,1]$ the restriction $T \bigl|_{L_2(A)}$ is not an isomorphic embedding. Then $T$ is not a sign embedding in the sense of Definition~\ref{se}.
\end{prop}

First, we note  that Theorem~\ref{L2main} is an immediate  corollary of Proposition~\ref{2notsignpr}, Lemma~\ref{signsomewhatnarrow}, and Theorem~\ref{A}.

\begin{proof}[Proof of Proposition~\ref{2notsignpr}]
Suppose that  there exists an operator $T:L_2\to L_2$  which  satisfies assumptions of Proposition~\ref{2notsignpr} and is a sign embedding in the sense of Definition~\ref{se}, i.e. there exists $\de>0$ so that for every sign $x\in L_2$
$$\|T(x)\|_2\ge 2\de \|x\|_2.$$

By Proposition~\ref{structure}, for $\e=\de/4$,  $X = L_2$ and $x_n = e_n$, a fixed orthonormal basis of $L_2$, there exist an operator $S:L_2\to L_2$ and a measure preserving homeomorphism $\sigma:[0,1]\to[0,1]$, which satisfy all conditions of Proposition~\ref{structure}, in particular $|a_n|\ge \de$ for all $n\ge 2$. Since every normalized block basis of $(e_n)$ is isometrically equivalent to $(e_n)$ itself (see e.g. \cite[Proposition~2.a.1]{LTI}), we may and do assume that $u_n = e_n$.

Let $B=\sigma^{-1}([0,\frac12))$. Since $T\bigl|_{L_2(B)}$ is not an isomorphic embedding, there exists $x\in L_2([0,\frac12))$ such that $\|x\|=1$ and $\|T(x\circ\sigma)\|\le \frac{\de}4$. By (3) of Proposition~\ref{structure}, we get that
\begin{equation}\lb{normSx}
\|Sx\|\le \frac\de2.
\end{equation}

Recall that the $L_2$-normalized Haar system $(h_n)_{n\in\bbN}$ is an orthonormal basis for $L_2[0,1]$ (see e.g. \cite[p.~128]{AK}). Thus there exist coefficients
$(b_n)_{n=0}^\infty$ such that
$$x=b_0{\bold{1}}_{[0,\frac12)}+\sum_{n\in M} b_nh_n,$$
where $M$ is the subset of $\bbN$ consisting of all $n$ such that $\supp h_n\subseteq[0,\frac12)$. Note that 
${\bold{1}}_{[0,\frac12)}=\frac12(h_1+h_2).$ Thus 
$$x=\frac12b_0h_1+\frac12b_0h_2 +\sum_{n\in M} b_nh_n.$$ 

Since $\|x\|=1$, and $(h_n)_{n\in\bbN}$ is an orthonormal basis, we get
\begin{equation}\lb{normx}
\frac14 b_0^2+\frac14 b_0^2+\sum_{n\in M} b_n^2=1.
\end{equation}

On the other hand, by \eqref{normSx} and properties of $S$ we have
\begin{equation*}
\begin{split}
\frac\de2\ge\|Sx\|&=\Big\|\frac12b_0a_2e_2 +\sum_{n\in M} b_na_ne_n\Big\|\\
&=\Big(\frac14 b_0^2a_2^2+\sum_{n\in M} b_n^2a_n^2\Big)^{\frac12}\\
&\ge \de\Big(\frac14 b_0^2+\sum_{n\in M} b_n^2\Big)^{\frac12}.
\end{split}
\end{equation*}

Thus
\begin{equation}\lb{small}
\frac14 b_0^2+\sum_{n\in M} b_n^2\le \frac14.
\end{equation}

Notice that \eqref{normx} and \eqref{small} imply that
$$\frac14 b_0^2\ge\frac34,$$
but \eqref{small} implies that 
$$\frac14 b_0^2\le\frac14.$$
The resulting contradiction ends the proof.
\end{proof}

\subsection{Final remarks}
It is easy to see that Problem~\ref{pr:1} has a negative answer for $2<p<\infty$. Indeed, this is witnessed by the operator $T$ defined in Example~\ref{Sch}. However, this operator is somewhat narrow for every $2<p<\infty$. Motivated by our proof of Theorem~\ref{L2main}, it is courious to ask whether for $2<p<\infty$ somewhat narrow operators could replace narrow operators in Problem~\ref{pr:1}, i.e. if the following problem is true.

\begin{prob} \label{pr:2}
Suppose $2 < p <\infty$, and an operator $T:L_p\to L_p$ is such that for every  measurable set $A \subseteq[0,1]$ the restriction $T \bigl|_{L_p(A)}$ is not an isomorphic embedding. Does it follow that $T$ is somewhat narrow?
\end{prob}

Note that for $1<p\le 2$, by Theorem~\ref{A}, Problems~\ref{pr:2} and \ref{pr:1} are identical.

We feel that an affirmative answer to Problem~\ref{pr:2} would be a very strong property of operators on $L_p$, $2 < p <\infty$, thus the author chose not to formulate Problem~\ref{pr:2} as a conjucture. A counterexample to Problem~\ref{pr:2} would also be very interesting.

\begin{ack}
This paper arose as a result of my participation in the  ``Problems and Recent Methods in Operator Theory Workshop'' held at the University of Memphis, October 15-16, 2015, with partial support by the National Science Foundation Award DMS-1546799.
I thank Professor Fernanda Botelho, the main organizer of the workshop, for the invitation, warm hospitality, and financial support, which allowed me to be a part of this very mathematically stimulating event.

I also thank Professor Mikhail Popov and the anonymous referee for valuable comments which improved the presentation of the paper.
\end{ack}


\begin{thebibliography}{12}
\bibitem{AK}
\newblock {\sc F. Albiac, N. J. Kalton}.
\newblock Topics in Banach space theory. 
\newblock Graduate Texts in Mathematics, 233. Springer, New York, 2006. xii+373 pp.

\bibitem{Bou81}
\newblock {\sc J.~Bourgain}.
\newblock {\em New classes of ${\mathcal L}_p$-spaces}.
\newblock Lect. Notes Math.
\newblock {\bf 889} (1981), 1--143.

\bibitem{BR83}
\newblock {\sc J. Bourgain, H. P. Rosenthal}.
\newblock {\em  Applications of the theory of semi-embeddings to Banach space theory}.
\newblock  J. Funct. Anal. 
\newblock {\bf 52, No 2} (1983), pp. 149--188.


\bibitem{DJS}
\newblock {\sc D.~Dosev, W.~B.~Johnson, G.~Schechtman}.
\newblock {\em Commutators on $L_p$, $1\le p<\infty$}.
\newblock J. Amer. Math. Soc.
\newblock {\bf 26,   No. 1} (2013), 101--127.

\bibitem{ES79}
\newblock {\sc P.~Enflo, T.~Starbird}.
\newblock {\em Subspaces of $L^1$ containing $L^1$}.
\newblock Studia Math.
\newblock {\bf 65, No 2} (1979), 203--225.

\bibitem{FR03}
\newblock {\sc J. Flores, C. Ruiz}.
\newblock {\em Domination by positive narrow operators}.
\newblock   Positivity 
\newblock {\bf 7} (2003), pp. 303--321.

\bibitem{GR83}
\newblock {\sc N. Ghoussoub, H. P. Rosenthal}.
\newblock {\em  Martingales, $G_\delta$-embeddings and quotients of $L^1$}.
\newblock    Math. Anal. 
\newblock {\bf 264, No 3} (1983), pp. 321--332.

\bibitem{JMST}
\newblock {\sc W.~B.~Johnson, B.~Maurey, G.~Schechtman, L.~Tzafriri}.
\newblock {\em Symmetric structures in Banach spaces}.
\newblock Memoirs of the Amer. Math. Soc.
\newblock {\bf 19, No 217} (1979).

\bibitem{KW76}
\newblock {\sc N. J. Kalton, A. Wilansky}.
\newblock {\em Tauberian operators on Banach spaces}. 
\newblock   Proc. Amer. Math. Soc. 
\newblock {\bf 57 no. 2} (1976),  251--255. 

\bibitem{K09}
\newblock {\sc I. V. Krasikova}.
\newblock {\em  A note on narrow operators in $L_\infty$}. 
\newblock   Math. Stud., 
\newblock {\bf 31, No 1} (2009), pp. 102--106.

\bibitem{LP71}
\newblock {\sc J. Lindenstrauss, A. Pe\l czy\'nski}.
\newblock {\em Contributions to the theory of the classical Banach spaces}. 
\newblock J. Functional Analysis 
\newblock {\bf 8 } (1971), 225--249. 

\bibitem{LTI}
\newblock {\sc J.~Lindenstrauss, L.~Tzafriri}.
\newblock { Classical Banach spaces, Vol. 1, Sequence spaces}.
\newblock Springer--Verlag, Berlin--Heidelberg--New York
\newblock (1977).

\bibitem{LPP79}
\newblock {\sc H. P. Lotz, N. T. Peck, H. Porta}.
\newblock {\em  Semi-embedding of Banach spaces}.
\newblock   Proc. Edinburgh Math. Soc. 
\newblock {\bf 22} (1979), pp. 233--240.

\bibitem{M1916}
\newblock {\sc D. Menchoff}. 
\newblock {\em Sur l'unicit\'e du d\'eveloppment trigonom\'etrique}.
\newblock Comptes Rendus de l'Acad\'emie des Sciences 
\newblock {\bf 163} (1916), 433--436. 


\bibitem{MP06}
{\sc V.~V.~Mykhaylyuk, M.~M.~Popov}.
\newblock {\em ``Weak" embeddings of $L_1$}.
\newblock Houston J. Math.
\newblock {\bf 32, No 4} (2006),  1139--1152.

\bibitem{MPRS}
\newblock {\sc V. Mykhaylyuk, M. Popov, B. Randrianantoanina, G. Schechtman}.
\newblock {\em Narrow and $\ell_2$-strictly singular operators from $L_p$}.
\newblock Israel J. Math. 
\newblock {\bf 203, no. 1} (2014) 81--108.


\bibitem{PP}
\newblock {\sc A.~M.~Plichko, M.~M.~Popov}.
\newblock {\em Symmetric function spaces on atomless probability spaces}.
\newblock Dissertationes Math. (Rozprawy Mat.)
\newblock {\bf 306} (1990), 1--85.



\bibitem{PRbook}
\newblock {\sc M.~M.~Popov, B. Randrianantoanina}.
\newblock  Narrow Operators on Function Spaces and Vector Lattices.
\newblock De Gruyter Studies in Mathematics 45,
\newblock Walter de Gruyter \& Co., Berlin, 2013. xiv+319 pp.




\bibitem{Ros81/82}
\newblock {\sc H.~P.~Rosenthal}.
\newblock {\em Some remarks concerning sign embeddings}.
\newblock Semin. D'analyse Fonct. Univ. of Paris VII.
\newblock (1981/82).

\bibitem{Ros83}
\newblock {\sc H.~P.~Rosenthal}.
\newblock {\em Sign-embeddings of $L^1$}.
\newblock Lect. Notes Math., Springer, Berlin
\newblock {\bf 995}  (1983), 155--165.

\bibitem{Ros84}
\newblock {\sc H.~P.~Rosenthal}.
\newblock {\em Embeddings of $L^1$ in $L^1$}.
\newblock Contemp. Math.
\newblock {\bf 26} (1984), 335--349.

\bibitem{T90}
\newblock {\sc M. Talagrand}.
\newblock {\em  The three space problem for $L^1$}.
\newblock J. Amer. Math. Soc.  
\newblock {\bf 3, No 1} (1990), pp. 9--29.



\end{thebibliography}
\end{document}